\documentclass{commat}

\usepackage{amsmath,amsthm,amsfonts,amscd,amssymb}

\theoremstyle{definition}

\theoremstyle{remark}

\title{Euclidean lattices: theory and applications}

\author{Lenny Fukshansky and Camilla Hollanti}

\authorinfo[L. Fukshansky]{Department of Mathematics, 850 Columbia Avenue, Claremont McKenna College, Claremont, CA 91711, USA}{lenny@cmc.edu}

\authorinfo[C. Hollanti]{Department of Mathematics and Systems Analysis, Aalto University, P.O. Box 11100, FI-00076 Aalto, Finland}{camilla.hollanti@aalto.fi}

\abstract{In this editorial survey we introduce the special issue of the journal Communications in Mathematics on the topic in the title of the article. Our main goal is to briefly outline some of the main aspects of this important area at the intersection of theory and applications, providing the context for the articles showcased in this special issue.}



\VOLUME{31}
\YEAR{2023}
\NUMBER{2}
\firstpage{251}
\DOI{https://doi.org/10.46298/cm.11596}

\begin{document}

\def\A{{\mathcal A}}
\def\AA{{\mathfrak A}}
\def\B{{\mathcal B}}
\def\C{{\mathcal C}}
\def\D{{\mathcal D}}
\def\F{{\mathcal F}}
\def\x{{\mathcal H}}
\def\I{{\mathcal I}}
\def\J{{\mathcal J}}
\def\K{{\mathcal K}}
\def\kk{{\mathfrak K}}
\def\L{{\mathcal L}}
\def\M{{\mathcal M}}
\def\O{{\mathcal O}}
\def\mm{{\mathfrak m}}
\def\MM{{\mathfrak M}}
\def\OO{{\mathfrak O}}
\def\R{{\mathcal R}}
\def\s{{\mathcal S}}
\def\V{{\mathcal V}}
\def\X{{\mathcal X}}
\def\Y{{\mathcal Y}}
\def\Z{{\mathcal Z}}
\def\H{{\mathcal H}}
\def\bee{{\mathbb B}}
\def\cee{{\mathbb C}}
\def\eee{{\mathbb E}}
\def\pee{{\mathbb P}}
\def\que{{\mathbb Q}}
\def\real{{\mathbb R}}
\def\zed{{\mathbb Z}}
\def\aaa{{\mathbb A}}
\def\ff{{\mathbb F}}
\def\kk{{\mathfrak K}}
\def\qbar{{\overline{\mathbb Q}}}
\def\kbar{{\overline{K}}}
\def\ybar{{\overline{Y}}}
\def\kkbar{{\overline{\mathfrak K}}}
\def\ubar{{\overline{U}}}
\def\eps{{\varepsilon}}
\def\ahat{{\hat \alpha}}
\def\bhat{{\hat \beta}}
\def\gt{{\tilde \gamma}}
\def\h{{\tfrac12}}
\def\ba{{\boldsymbol a}}
\def\bb{{\boldsymbol b}}
\def\be{{\boldsymbol e}}
\def\bei{{\boldsymbol e_i}}
\def\bc{{\boldsymbol c}}
\def\bm{{\boldsymbol m}}
\def\bk{{\boldsymbol k}}
\def\bi{{\boldsymbol i}}
\def\bl{{\boldsymbol l}}
\def\bq{{\boldsymbol q}}
\def\bu{{\boldsymbol u}}
\def\bt{{\boldsymbol t}}
\def\bs{{\boldsymbol s}}
\def\bv{{\boldsymbol v}}
\def\bw{{\boldsymbol w}}
\def\bx{{\boldsymbol x}}
\def\bX{{\boldsymbol X}}
\def\bz{{\boldsymbol z}}
\def\bwy{{\boldsymbol y}}
\def\bY{{\boldsymbol Y}}
\def\bL{{\boldsymbol L}}
\def\baa{{\boldsymbol\alpha}}
\def\bbb{{\boldsymbol\beta}}
\def\bet{{\boldsymbol\eta}}
\def\bxi{{\boldsymbol\xi}}
\def\bo{{\boldsymbol 0}}
\def\bol{{\boldkey 1}_L}
\def\ep{\varepsilon}
\def\p{\boldsymbol\varphi}
\def\q{\boldsymbol\psi}
\def\rank{\operatorname{rank}}
\def\aut{\operatorname{Aut}}
\def\lcm{\operatorname{lcm}}
\def\sgn{\operatorname{sgn}}
\def\spn{\operatorname{span}}
\def\md{\operatorname{mod}}
\def\Norm{\operatorname{Norm}}
\def\dim{\operatorname{dim}}
\def\det{\operatorname{det}}
\def\Vol{\operatorname{Vol}}
\def\rk{\operatorname{rk}}
\def\ord{\operatorname{ord}}
\def\ker{\operatorname{ker}}
\def\div{\operatorname{div}}
\def\Gal{\operatorname{Gal}}
\def\GL{\operatorname{GL}}
\def\Sym{\operatorname{Sym}}

\tableofcontents

\section{Introduction}
\label{intro}

A {\it lattice} $L$ in a Euclidean $n$-dimensional space $\eee_n$ is a discrete subgroup of rank $1 \leq k \leq n$. This is equivalent to saying that there exists a collection of linearly independent elements $\ba_1,\dots,\ba_k \in \eee_n$ (always written as column vectors) such that
$$L = \left\{ \sum_{i=1}^k c_i \ba_i : c_1,\dots,c_k \in \zed \right\} = A\zed^k,$$
where $\ba_1,\dots,\ba_k$ is a basis for $L$ and $A = (\ba_1\ \dots\ \ba_k)$ is the corresponding $n \times k$ basis matrix. If this is the case, then for any $U \in \GL_k(\zed)$,
$$L = A\zed^k = (AU)\zed^k,$$
and so $AU$ is again a basis matrix for $L$. Identifying $\eee_n$ with the real space $\real^n$, we can therefore identify the space of all rank-$k$ lattices in $\real^n$ with the space $\GL_k(\real)/\GL_k(\zed)$ of all orbits of $\GL_k(\real)$ under the action of $\GL_k(\zed)$ by right multiplication.
\smallskip

Theory of Euclidean lattices connects number theory to convex and discrete geometry. The study of lattices originally emerged as an important subject in connection with classical discrete optimization problems like sphere packing, covering and kissing number problems, dating as far back as the celebrated 1611 conjecture of Kepler and even earlier; see the classical books of Conway \& Sloane \cite{conway} and of Martinet \cite{martinet} for a fairly comprehensive exposition of lattice theory and its many connections, as well as \cite{szpiro} for a popular account of the fascinating history of Kepler's conjecture. Lattices have really come into their own in the context of Minkowski's geometry of numbers (see \cite{minkowski} for Minkowski's original treatise, as well as the standard books \cite{cass:geom} by Cassels, \cite{lek} by Gruber \& Lekkerkerker and \cite{gruber} by Gruber for more contemporary accounts).
\smallskip

Theory of lattices has seen some very exciting developments and applications over the last century, including Minkowski's proof of the finiteness of class number, major results in the arithmetic theory of quadratic forms, advances in discrete and convex geometry and optimization, Diophantine approximations, geometric combinatorics, coding theory, cryptography, and many other areas of mathematics. The recent decades have, in particular, seen such major breakthroughs as the proof of Kepler's conjecture by Hales \& Ferguson \cite{hales}, affirming that the face-centered-cubic lattice provides the densest sphere packing in dimension 3, as well as the spectacular results by Viazovska \emph{et al.} \cite{viaz1}, \cite{viaz2} on the optimality of $E_8$ and the Leech lattice for packing density in dimensions 8 and 24, respectively (Maryna Viazovska received a Fields medal for this work in 2022).
\smallskip

The main goal of our special issue is to collect in one place several of the recent developments and expository surveys on the various aspects of lattice theory and its applications. In the following sections, we will briefly introduce a few different facets of this theory and indicate how different contributions of this special issue fit into the general framework.
\bigskip

\section{Geometry of numbers and Diophantine approximations}
\label{gn}

The first essential invariant of the lattice $L$ as above is its {\it determinant}, which is defined as
$$\det(L) := \sqrt{\det(A^{\top} A)}$$
for any choice of a basis matrix $A$: this is well-defined, since $|\det(U)| = 1$ for any $U \in \GL_k(\zed)$. Analytically, this is the volume of a fundamental parallelotope
$$\left\{ A\bx : \bx \in [0,1)^k \right\},$$
which is a full set of coset representatives for the quotient group $V/L$, where $V = \spn_{\real} L$. In fact, $\det(L)$ is the volume of the closure of any such fundamental domain, including the important {\it Voronoi cell}
\[\V(L) = \left\{ \bx \in V : \|\bx\| \leq \|\bx-\bwy\|\ \forall\ \bwy \in L \right\},\]
i.e., the set of all points in $V$ that are no further from the origin than from any other point of the lattice. Our lattice $L$ has full rank in the $k$-dimensional subspace $V$ of $\real^n$, which can be identified with the Euclidean space $\real^k$. As such, we will only talk of full-rank lattices in $\real^n$ from now on.
\smallskip

We can now define the {\it sphere packing} and the {\it sphere covering} associated to $L$: inscribe a closed ball $B_1$ in $\real^n$ of maximal possible radius into $\V(L)$ and circumscribe a ball $B_2$ of minimal possible radius around $\V(L)$, then translating $\V(L)$ by all the points of $L$ we obtain a packing of non-overlapping translates of $B_1$ in $\real^n$ and a covering of $\real^n$ by translates of $B_2$. Hence the radius of $B_1$ is called the {\it packing radius} $r(L)$ of $L$ and the radius of $B_2$ the {\it covering radius} $R(L)$ of~$L$. Now, the {\it packing density} $\delta(L)$ and the {\it covering thickness} $\Theta(L)$ are given by the formulas
\[\delta(L) = \frac{\Vol_n(B_1)}{\Vol_n(\V(L))} = \frac{\omega_n r(L)^n}{\det(L)},\ \Theta(L) = \frac{\Vol_n(B_2)}{\Vol_n(\V(L))} = \frac{\omega_n R(L)^n}{\det(L)},\]
where $\omega_n$ is the volume of a unit ball $\bee_n$ in $\real^n$. In fact, these radii are closely related to another collection of important invariants of the lattice $L$, called successive minima. 

Let $K$ be a closed convex $\bo$-symmetric set of positive volume in $\real^n$. The {\it successive minima} of the lattice $L$ with respect to $K$,
$$0 < \lambda_1(L,K) \leq \dots \leq \lambda_n(L,K),$$
are defined as
\[\lambda_i(L,K) = \min \left\{ t \in \real_{>0} : \dim_{\real} \spn_{\real} L \cap tK \geq i \right\},\]
i.e., the smallest real number $t$ so that the homogeneous expansion of $K$ by a factor of $t$ contains at least $i$ linearly independent points of $L$. In the special case when $K$ is the unit ball $\bee_n$ centered at the origin in $\real^n$, we refer to $\lambda_i(L,K)$ simply as $\lambda_i(L)$, the successive minima of the lattice. It is then not difficult to see that the packing radius is precisely half the distance from the origin to the shortest nonzero lattice point, i.e.
$$r(L) = \frac{1}{2} \lambda_1(L).$$
The more delicate inequalities of Jarnik (see, e.g., \cite[Section 13.2, Theorem 1 and Theorem 4]{lek}) also assert that the covering radius satisfies
$$\frac{1}{2} \lambda_n(L) \leq R(L) \leq \frac{1}{2} \sum_{i=1}^n \lambda_i(L).$$

Successive minima have been studied extensively by Minkowski himself and by many other mathematicians working in number theory, discrete and convex geometry, and even analysis. In particular, Minkowski's inequalities on successive minima state that
\begin{equation}
\label{mink2}
\frac{2^n \det(L)}{n! \Vol_n(K)} \leq \prod_{i=1}^n \lambda_i(L,K) \leq \frac{2^n \det(L)}{\Vol_n(K)}.
\end{equation}
The survey paper \cite{martin_iskander} by I. Aliev \& M. Henk in this special issue gives an overview of the impact of successive minima on convex and discrete geometry. One significant application of successive minima inequalities that the authors discuss is Siegel's lemma, a vital tool in Diophantine approximations and transcendental number theory, which provides a bound on the size of a ``smallest" nonzero solution (or, more generally, a collection of such solutions) to a system of linear forms over a given ring or field of arithmetic interest. The fact that such a solution exists over a given field is guaranteed by the assumption that these linear forms are linearly dependent over the same field.
\smallskip

On the other hand, assume that we have a system of linear forms that are linearly independent over $\que$. Then they will not be simultaneously equal to zero at any nonzero point of the integer lattice. A natural question in Diophantine approximations is how small can such a collection of linear forms in $n$ variables simultaneously be on $\zed^n \setminus \{\bo\}$? This question can be made precise by studying the extreme values of certain appropriately defined exponents of approximation, which is done in the paper~\cite{german} by O. German in our special issue. His main tool is a lemma of Davenport on successive minima. At the end of this paper, a question about the spectra of these newly introduced Diophantine exponents is formulated.
\bigskip

\section{Special classes of lattices}
\label{scl}

As we remarked above, the analogues of Kepler's conjecture on the densest possible sphere packing in dimension $3$ has also been proved in dimensions $8$ and $24$. In fact, the optimal sphere packing has been obtained earlier in dimension $2$ by L. Fejes T\'oth \cite{toth}, who gave the first complete proof of what was previously known as Thue's theorem. These are all the dimensions (besides the trivial dimension $1$) in which the optimal sphere packings are known. If, however, we restrict our consideration to lattice packings only, then the optimal results are known in dimensions $1 \leq n \leq 8$ and $n=24$ (see \cite{conway}).

One can then pose a natural question: what properties should a lattice possess to be a potential candidate for maximizing lattice packing density in its dimension? From our discussion above, it is evident that the packing density of a full-rank lattice $L$ in $\real^n$ is given by the formula
$$\delta(L) = \frac{\omega_n \lambda_1(L)^n}{2^n \det(L)}.$$
Let us define an equivalence relation of similarity on lattices in $\real^n$ as follows: two lattices $L_1$ and $L_2$ are called {\it similar} if there exists a positive real constant $\alpha$ and an $n \times n$ real orthogonal matrix $U$ so that $L_2 = \alpha U L_1$. In this case, it is easy to see that $\delta(L_1) = \delta(L_2)$, and so the packing density is constant on a given similarity class. Hence, restricting to unimodular lattices (determinant $=1$) we can write
$$\delta(L) = \frac{\omega_n}{2^n}\ \lambda_1(L)^n,$$
meaning that maximizing packing density is equivalent to maximizing the first successive minimum $\lambda_1(L)$. By Minkowski's inequalities~\eqref{mink2}, the product of all successive minima in this case is bounded by dimensional constants:
$$\frac{2^n}{n!\ \omega_n} \leq \prod_{i=1}^n \lambda_i(L) \leq \frac{2^n}{\omega_n},$$
where $0 < \lambda_1(L) \leq \dots \leq \lambda_n(L)$. Thus the first step in the direction of maximizing $\lambda_1(L)$ is to take a lattice $L$ with all successive minima equal: lattices like this are called {\it well-rounded (WR)}. This property is preserved under similarity, so we can talk of WR similarity classes, of which there are infinitely many in any dimension $n \geq 2$. There has been quite a bit of work in recent years on various explicit algebraic constructions of WR lattices. Some most notable such constructions come from ideals in algebraic number fields via Minkowski embedding into Euclidean space, the so-called {\it ideal lattices}. As such, an interesting question remains: under which conditions does an ideal in a number field give rise to WR ideal lattices? A detailed study of WR ideal lattices has been initiated in~\cite{lf_kp}. While this question has been answered for quadratic number fields and for some special families of number fields of higher degree, in general, it is wide open. In article~\cite{ha_tran} in our special issue, D.T. Tran, N. H. Le, and H. T. N. Tran conduct a thorough investigation and establish conditions for the existence of WR ideal lattices coming from cyclic number fields of degree $3$ and $4$. Their paper starts out with a brief overview of the previous results in this area and also contains a fairly extensive bibliography. 
\smallskip

The packing density function is continuous on $\GL_n(\real)/\GL_n(\zed)$, the space of full-rank lattices in $\real^n$, and hence we can talk about its local extrema on this space. While the WR condition is necessary for a local maximum to be achieved, this condition is not sufficient. Define the set {\it of minimal vectors} of a lattice $L$ as
$$S(L) = \left\{ \bx \in L : \|\bx\| = \lambda_1(L) \right\},$$
and let $m$ be the cardinality of $S(L)$. Notice that $m$ is necessarily even since minimal vectors come in $\pm$ pairs (more generally, $m$ is divisible by the order of the group of linear automorphisms of $L$ since it acts on $S(L)$ by left multiplication). Further, if $L$ is WR then $m \geq 2n$. Lattice $L$ is called {\it eutactic} if there exist positive real coefficients $c_1,\dots,c_m$ such that for any vector $\bv \in \real^n$,
$$\|\bv\| = \sum_{i=1}^m c_i \left( \bv^{\top} \bx_i \right)^2,$$
where $S(L) = \{ \bx_1,\dots,\bx_m \}$. On the other hand, a lattice $L$ is called {\it perfect} if the space of $n \times n$ real symmetric matrices $\Sym_n(\real)$ can be spanned (as a real vector space) by symmetric matrices coming from the minimal vectors of $L$, i.e.
$$\Sym_n(\real) = \spn_{\real} \left\{ \bx \bx^{\top} : \bx \in S(L) \right\}.$$
Since $\dim_{\real} \Sym_n(\real) = \frac{n(n+1)}{2}$ and for any vector $\bx \in S(L)$, $\bx \bx^{\top} = (-\bx) (-\bx)^{\top}$, this perfection condition implies that the cardinality $m$ of $S(L)$ is at least $n(n+1)$. The eutaxy and perfection conditions on lattices are independent (i.e., there are eutactic non-perfect lattices and there are perfect non-eutactic lattices) and they are both preserved under similarity. Furthermore, there are only finitely many eutactic and finitely many perfect similarity classes in any given dimension, although their number grows very fast with the dimension (for instance, for sufficiently large $n$ the number of perfect similarity classes in $\real^n$ is $> e^{n^{1-\eps}}$ for any $\eps > 0$; see~\cite{bacher}). Both, perfect and eutactic lattices are necessarily WR and a famous theorem of Voronoi (1908) asserts that a lattice corresponds to a local maximum of the packing density function in its dimension (called {\it extreme} lattice) if and only if it is perfect and eutactic (see, e.g., \cite{martinet}). This observation drives an interest in the classification of perfect lattices, a subject of active research often pursued in the language of quadratic forms.

For $L = A\zed^n$, the Euclidean norm of any vector $\bx = A\bwy \in L$ can be computed~as
$$\|\bx\|^2 = Q_A(\bwy) := \bwy^{\top} \left( A^{\top} A \right) \bwy,$$
where $Q_A(\bwy)$ is a positive definite quadratic form with a symmetric coefficient matrix $A^{\top} A$. Quadratic forms corresponding to different bases of the same lattice are called {\it arithmetically equivalent}: they have the same spectrum of values on $\zed^n$. There is a bijective correspondence between positive definite arithmetic equivalence classes of quadratic forms in $n$ variables and lattices in $\real^n$. The symmetric coefficient matrix of a positive definite quadratic form is then called perfect if the corresponding lattice is perfect. The paper by V. Dannenberg and A. Sch\"urmann \cite{achill} in our special issue builds on the classical theory of such perfect matrices to introduce and initiate a study of their generalization, the so-called perfect copositive matrices: a matrix $B \in \Sym_n(\real)$ is called {\it copositive} if
$$\bwy^{\top} B \bwy \geq \bo$$
for all $\bwy$ in the positive orthant $\real^n_{\geq 0}$ (in contrast to the usual positive definite matrices satisfying $\bwy^{\top} B \bwy \geq \bo$ for all nonzero $\bwy \in \real^n$). The authors look at the cone of copositive matrices and study the distribution of perfect matrices in this cone.

\bigskip

\section{Arithmetic of quadratic forms}
\label{quad}

As indicated above, the study of lattices is intrinsically connected to the arithmetic theory of quadratic forms. A key question in that theory is that of representation. A quadratic form $Q(\bwy)$ in $n$ variables can always be written in the form
$$Q(\bwy) = \bwy^{\top} B \bwy,$$
where $B$ is a real symmetric coefficient matrix. This form $Q$ is called {\it integral} if $Q(\bwy) \in \zed$ for every $\bwy \in \zed^n$ and it is called {\it classically integral} if $B$ is an integer matrix; notice that this second property is stronger than the first. An integral form $Q$ is said to {\it represent} an integer $m$ if there exists $\bwy \in \zed^n$ such that $Q(\bwy) = m$, and $Q$ is said to be {\it universal} if it represents every positive integer $m$. This is equivalent to the corresponding lattice containing vectors of every possible (squared) integer Euclidean norm. 

Perhaps the starting point of the theory of universal quadratic forms is the famous classical theorem of Lagrange (1770) stating that the positive definite integral quadratic form given by the sum of four squares is universal (see, for instance, \cite{weil} for details). On the other hand, no positive definite integral form in fewer than four variables can be universal. The major results on universal forms from the past thirty years include the impressive necessary and sufficient universality criteria for integral and classically integral quadratic forms in any number of variables, known as theorems 290 and 15, respectively. 

The survey article by V. Kala \cite{vita} (based on the author's lectures on this subject) in this special issue gives an overview of the theory of universal quadratic forms, including these celebrated theorems, but placing the main focus on the recent developments for quadratic lattices over the ring of integers $\O_K$ in a totally real number field $K$. The key tool emphasized by the author is the notion of an indecomposable element in $\O_K$: a totally positive element in $\O_K$ is called {\it indecomposable} if it cannot be written as a sum of two other totally positive elements in the same ring. The significance of indecomposables in the context of quadratic forms is that they essentially appear as coefficients of diagonal universal forms over $\O_K$ and hence the number of their square classes gives a lower bound on the number of variables in which such forms can exist. The author carefully develops the theory of indecomposables in this context, showing also some interesting connections, including one to continued fractions.

\bigskip

\section{Geometric combinatorics and integer geometry}
\label{gcomb}

Another important facet of the theory has to do with counting lattice points in compact domains in $\real^n$. More specifically, let us start with a full-rank lattice $L$ and a compact measurable set $K \subset \real^n$ of positive volume. Let $r \in \real_{>0}$ and define the counting function 
$$f_{L,K}(r) = \left| L \cap rK \right|.$$
One can ask how does $f_{L,K}(r)$ grow as $r \to \infty$? The first observation is that each point $\bx \in L$ is contained in its unique translate of the Voronoi cell $\bx+\V(L)$, hence counting lattice points can be replaced by counting translates of the Voronoi cell. As $r$ becomes large, the number of such translates that are fully contained in $rK$ gives the main term of the asymptotic formula for $f_{L,K}(r)$, whereas the error term comes from the number of such translates intersecting tie boundary of $rK$ whose corresponding lattice points are in $rK$. Hence the main term can be approximated simply by the quotient of the volume of $rK$ by the volume of the Voronoi cell, $\det(L)$. Under appropriate smooth conditions on the boundary of $K$, such as Lipschitz parametrizability, the error term can be controlled, and the following asymptotic holds (see, e.g., \cite{lang}, Chapter VI, \S 2, Theorem 2):
$$f_{L,K}(r) = \frac{\Vol_n(K)}{\det(L)}\ r^n + O(r^{n-1}).$$
A considerably more delicate problem is to give tight (and as explicit as possible) estimates on the error term
$$\left| f_{L,K}(r) - \frac{\Vol_n(K)}{\det(L)}\ r^n \right|.$$
There is a vast amount of literature on different versions of this counting problem. In fact, this problem is not fully resolved even in a seemingly simple case of $L=\zed^2$ and $K$ being the unit circle $S^1$ -- this is the famous Gauss circle problem, where the standing conjecture is that
$$\left| f_{\zed^2,S^1}(r) - \pi r^2 \right| = O(r^{1/2+\eps})$$
for any $\eps > 0$. 

A variation of this counting problem is treated in a paper of J. D. Vaaler~\cite{vaaler} in this special issue: given an $n \times m$ real matrix $A$, $m \leq n$, obtain an estimate on the error term
$$\left| f_{A\zed^m,rB(\bx)} - \Vol_m(B(\bx))\ r^m \right|$$
for the number of points of the lattice $A\zed^m$ in the ball $rB(\bx)$ of radius $r$ centered at an arbitrary point $\bx$ in the subspace $A\real^m \subseteq \real^n$ spanned by this lattice, as $r \to \infty$. While a number of estimates on such quantities have been previously obtained (see \cite{vaaler} for some bibliography), the author's estimate is explicit and uniform over all matrices $A$ with norm bounded by an explicit constant. Further, his inequality takes a particularly simple form for dimension $m=3$. The author's method uses careful analysis of extremal functions; as such, Bessel functions naturally occur in the estimates.
\smallskip

The situation with counting integer lattice points becomes more manageable when the compact set $K$ is a convex lattice polytope. Indeed, assume $K$ is a convex polytope in $\real^n$ with positive volume and vertices at points of the integer lattice $\zed^n$. Consider the counting function $f_{\zed^n,K}(r)$ for integer values of the homogeneous expansion parameter $r$. A classical theorem of Ehrhart (1962) states that $f_{\zed^n,K}(r)$ is a polynomial in $r$ of degree $n$ with integer coefficients, where the leading coefficient is $\Vol_n(K)$ (see \cite{beck_robins} as well as Chapter 12 of \cite{miller} for a nice introduction to Ehrhart theory). This polynomial is called {\it Ehrhart polynomial} of the polytope $K$. More generally, let us define the {\it integer point transform} of the polytope $rK$ by
$$\sigma_{rK}(\bxi) = \sum_{\bv \in \zed^n \cap rK} e^{2\pi i (\bv^{\top} \bxi)},$$
for all $\bxi \in \real^n$. In particular, notice that 
$$\sigma_{rK}(\bo) = \sum_{\bv \in \zed^n \cap rK} 1 = f_{\zed^n,K}(r),$$
and hence the integer point transform of a polytope is a certain generalization of Ehrhart polynomial. In his paper \cite{sinai} in this special issue, S. Robins proves that the integer point transform is a complete invariant of the polytope in the following sense: two lattice polytopes $K_1$ and $K_2$ are equal to each other if and only if $\sigma_{K_1}(\bxi^*) = \sigma_{K_2}(\bxi^*)$ for
$$\bxi^* = \frac{1}{\pi} (\sqrt{p_1},\dots,\sqrt{p_n})^{\top},$$
where $p_1,\dots,p_n$ are the first $n$ primes. In fact, the author first uses the classical Lindemann–Weierstrass theorem
from transcendental number theory to prove the analogous property for equality of arbitrary finite sets of integer lattice points instead of sets contained in polytopes and then passes to polytopes. Further, he proves the complete invariant property also for Fourier transforms of general rational polyhedra. Additionally, he discusses lattice-spanning properties of polytopes and the integer point transform of finite abelian groups.
\smallskip

Ehrhart's theorem is often seen as a higher-dimensional generalization of the famous Pick's theorem (1899; see, e.g., \cite{beck_robins} and Chapter 2 of \cite{karpenkov_book}): if $S$ is the area of an integer polygon in the plane, $I$ is the number of integer lattice points in its interior and $E$ is the number of integer lattice points on its boundary, then
$$S = I + \frac{E}{2} - 1.$$
The essential feature of this theorem is that it connects a combinatorial notion (the number of integer lattice points in a polygon) with an analytic notion (the area of this polygon). This is the main idea of integer geometry: introducing ``discrete" ways of measuring some traditionally ``continuous" objects, as alluded to in the title of Beck \& Robins's book \cite{beck_robins}. A good introduction to integer geometry and its connections to continued fractions is given in Karpenkov's book \cite{karpenkov_book}. In their paper \cite{karpenkov} in this special issue, J. Blackman, J. Dolan and O. Karpenkov take this exploration a step further and introduce the theory of multidimensional integer trigonometry. The integer length and integer area are defined in terms of indices of sublattices in the integer lattice generated by lattice points on a given line segment or in a given triangle (which then generalizes to arbitrary polygons via sums over triangulations). Integer area can then be used to define integer trigonometric functions in the plane, which are also closely connected to continued fractions. After giving a careful exposition of planar integer trigonometry, the authors of \cite{karpenkov} present a generalization of this theory to higher dimensions via integer volume of appropriate simplices. They prove a variety of different properties of integer trigonometric functions in arbitrary dimensions and discuss an algorithmic approach to constructions of rational polyhedra via given collections of rational cones. They also discuss approximations of simplicial cones, which generalize classical approximation by continued fractions.

\bigskip

\section{Applications to coding theory and cryptography}
\label{appl}

Arithmetic lattices have also made their way into many applications, perhaps most notably within coding theory and cryptography. 

\subsection{Lattices from error-correcting codes}
The association of lattices with error-correcting codes is natural and, in order to reduce the decoding complexity, a possible direction is the construction of multilevel lattices from a family of nested codes, allowing for \emph{multistage decoding}. Several different constructions have been used to derive lattices from codes \cite{conway}. To provide one explicit example, let $\rho:\mathbb{Z}_q\rightarrow \mathbb{Z}$ be the standard inclusion map, which can be naturally extended to vectors and matrices. Then, the $q$-ary Construction A lattice associated to the linear code $C\subseteq \mathbb{Z}_q^n$ can be defined as 
$$
L_A(C)=\rho(C)+q\mathbb{Z}^n.
$$

In the article \cite{costa} in this special issue, F. do Carmo Silva, A. P. de Souza, E. Strey, and S. I. R. Costa consider Constructions D, D$^\prime$, and A from nested $q$-ary linear codes over $\mathbb{Z}_q$. They study the volume,
$L_P$-minimum distance ($1 \leq P \leq\infty$), and lower bounds for the coding gain of these constructions. Further, the aforementioned multistage decoding method is extended with re-encoding to Construction D$^\prime$ from $q$-ary linear codes under specific conditions. The definitions of Constructions D and D$^\prime$ are somewhat more involved, and we refer the reader to the article for more details.

\subsection{Lattice-based cryptography}
One of the most promising paradigms for post-quantum security is \emph{lattice-based cryptography}, often based on different variants of the so-called \emph{learning with errors (LWE)} problem \cite{regev}. The hardness of such cryptosystems can be proved by providing a reduction from a known hard lattice problem, e.g., the approximate shortest vector problem. 

To give an example, let us consider the ring $R_q=\mathbb{F}_q[x]/(f(x))$, where $q$ is prime and $f(x)\in\mathbb{Z}[x]$ is monic and irreducible. The polynomial learning with errors (PLWE) \cite{rosca} decision problem asks to distinguish, with a non-negligible advantage, a sample $(a,b=as+e)\in R_q^2$, where $s$ and $e$ are drawn from an appropriate discrete Gaussian distribution, from a uniformly random sample $(a,b)\in R_q^2$. In article \cite{ivan} in our special issue,  I. Blanco-Chac\'on, R. Dur\'an-Di\'az, R. Njah Nchiwo, and B. Barbero-Lucas study a decisional attack against a version of the PLWE problem in which the samples are taken from a certain proper subring of large dimension of the cyclotomic ring $\mathbb{F}_q[x]/(\Phi_{p^k}(x))$ for $k>1$, $\Phi$ not totally split, and $q\equiv 1$ (mod $p$). The attack exploits the fact that the roots of $\Phi$ have zero traces over suitable sub-extensions. This allows for a good attack success probability as a function of input samples. The paper points out a nice open question regarding the existence of rings with the related distribution-respecting reduction map. We refer to the article for more details as well as for an exposition on the ring and polynomial LWE problems.

\subsection{Lattice codes for secure wireless communications}
Yet another interesting application of arithmetic lattices appears in the context of \emph{physical layer security}. Namely, lattice coset codes can be utilized for communication over the wireless medium, where eavesdroppers may receive the transmitted signals in addition to the legitimate receiver \cite{oggier}. The security of such physical layer communications can be measured in many different ways, including the \emph{eavesdropper's correct decision probability} or the \emph{information leakage}. It has been shown that both of these quantities are bounded from above by the so-called \emph{flatness factor} \cite{belfiore}, yielding a natural criterion for the flatness factor of the lattice to be minimized. Essentially, the flatness factor $\epsilon_L(\sigma)$ measures the deviation of the lattice Gaussian distribution from the uniform distribution on a Voronoi cell, and it is closely related to the lattice theta series
$
\Theta_L(q)=\sum_{x\in L}q^{||x||^2}
$
as follows:
\[\epsilon_L(\sigma)=\frac{\Vol(L)}{(\sqrt{2}\pi\sigma)^n}\Theta_L(e^{-\frac{1}{2\sigma^2}})-1=\Theta_{L^*}(e^{-2\pi\sigma^2})-1,\]
where $L^*$ denotes the dual lattice. 

 For a ``flat'' lattice, it is harder to distinguish the received message from a uniformly random sample. In order to minimize the flatness factor, well-rounded lattices have been proposed as a coding solution \cite{cami}. This motivates the search for good well-rounded lattices in small and moderate dimensions, in addition to the purely theoretical interest. In this special issue, well-rounded ideal lattices from cyclic cubic and quartic fields are studied in article \cite{ha_tran}, as already mentioned in Section~\ref{scl}.
\bigskip

{\bf Acknowledgement:} We wish to thank the referees for the thorough reading and helpful comments. Fukshansky was partially supported by the Simons Foundation grant \#519058. Hollanti was partially supported by the Academy of Finland grant \#351271 and by the Finnish Ministry of Defence MATINE grant \#2500M-0147.
\bigskip

{\small\bibliography{commat}}

\EditInfo{July 18, 2023}{August 26, 2023}{Ivan Kaygorodov}
\end{document}